\documentclass[a4paper]{amsart}
\usepackage{amsfonts}
\theoremstyle{plain}
\begingroup
\newtheorem{theorem}{Theorem}[section]
\newtheorem{lemma}[theorem]{Lemma}
\newtheorem{proposition}[theorem]{Proposition}
\newtheorem{corollary}[theorem]{Corollary}
\newtheorem{remark}[theorem]{Remark}
\newtheorem{definition}[theorem]{Definition}
\theoremstyle{definition}
\theoremstyle{remark}
\numberwithin{equation}{section}

\newcommand{\hs}{{\mathcal H}}
\newcommand{\ks}{{\mathcal K}}

\newcommand{\leb}{{\mathcal L}}
\newcommand{\ds}{{\mathcal D}}

\newcommand{\Ps}{{\mathcal P}}

\newcommand{\R}{{\mathbb R}}
\newcommand{\N}{{\mathbb N}}

\newcommand{\msim}{{\rm M}^{n\times n}_{\rm sym}}

\newcommand{\Om}{\Omega}
\newcommand{\Omb}{\overline{\Omega}}

\newcommand{\diam}[1]{{\rm diam}(#1)}
\newcommand{\supt}[1]{{\rm supp}(#1)}
\newcommand{\tint}[1]{{\rm int}(#1)}
\newcommand{\tintpi}[1]{{\rm int}_{\pi}(#1)}

\title
[A STABILITY RESULT FOR NEUMANN PROBLEMS IN DIMENSION $N \ge 3$]
{A STABILITY RESULT FOR NEUMANN PROBLEMS \\
IN DIMENSION $N \ge 3$}
\author[A. Giacomini]
{Alessandro Giacomini}
\address[Alessandro Giacomini]{S.I.S.S.A., Via Beirut 2-4, 34014, Trieste,
Italy}
\email[A. Giacomini]{giacomin@sissa.it}
\begin{document}
\vskip .2truecm
\begin{abstract}
\small{We give a sufficient condition in dimension $N \ge 3$ in order to obtain
the stability of a sequence of Neumann problems on fractured domains.}
\end{abstract}
\maketitle

\section{Introduction}
\label{intr}

Given $\Om$ open and bounded in $\R^N$, $(K_n)$ a sequence of compact sets
in $\R^N$, consider the following Neumann problems
\begin{equation}
\label{intpbn}
\left\{
\begin{array}{l}
-\Delta u + u = f \\ \frac{\partial u}{\partial \nu}=0
\end{array}
\begin{array}{l}
{\rm in}\; \Om \setminus K_n \\ {\rm on}\; \partial \Om \cup
(\partial K_n \cap \Om)
\end{array}
\right.
\end{equation}
with $f \in L^2(\Om)$: we intend (\ref{intpbn}) satisfied in the usual weak
sense of Sobolev spaces, that is $u \in H^1(\Om \setminus K_n)$ and
$$
\int_{\Om \setminus K_n} \nabla u \nabla \varphi
+\int_{\Om \setminus K_n} u\varphi=
\int_{\Om \setminus K_n}f\varphi
$$
for all $\varphi \in H^1(\Om \setminus K_n)$.
If $(K_n)$ converges to a compact set $K$ in the Hausdorff metric , we look for
conditions on the sequence $(K_n)$ such that, considered the problem
\begin{equation}
\label{intpbnlim}
\left\{
\begin{array}{l}
-\Delta u + u = f \\ \frac{\partial u}{\partial \nu}=0
\end{array}
\begin{array}{l}
{\rm in}\; \Om \setminus K \\ {\rm on}\; \partial \Om \cup
(\partial K \cap \Om),
\end{array}
\right.
\end{equation}
the solutions $u_n$ of (\ref{intpbn}) (extended to $0$ on $K_n \cap \Om$)
converge to the solution $u$ of (\ref{intpbnlim}) (extended to $0$ on
$K \cap \Om$). If this is the case, we say that the Neumann problems
(\ref{intpbn}) are stable.
\par
The problem of stability for elliptic problems
under Neumann boundary conditions has been widely investigated.
Usually, since in general the domains
$\Om \setminus K_n$ are not regular, it is not possible to deal with the problem
using extension operators (see for example \cite{Kru2}, \cite{MarKru}).
\par
In dimension $N=2$, Chambolle
and Doveri \cite{ChD} in 1997 proved a stability result under a uniform
limitation of $\hs^1(K_n)$ and of the number of the connected components of $K_n$;
Bucur and Varchon \cite{Buc} in 2000 proved that
if $K_n$ has at most $m$ connected components ($m \in \N$), the stability of the
problems is equivalent to the condition
$\leb^2(\Om \setminus K_n) \to \leb^2(\Om \setminus K)$.
\par
In dimension $N \ge 3$, the bound on the number of the connected components of
$K_n$ is not a relevant feature and a condition similar to that  of Bucur and
Varchon
doesn't hold: in fact, problems (\ref{intpbn}) could be not stable even if
the sets $K_n$ are connected. In 1997, Cortesani \cite{Co} proved that in general,
if $K$ is contained in a $C^1$ submanifold of $\R^N$, the limit of solutions of
(\ref{intpbn}) satisfies a transmission condition on $K$.
Several results on this transmission condition are known under additional
assumptions on $(K_n)$. In the case in which $K_n$ is contained in a hyperplane
$M$ and is the complement in $M$ of a periodic grid of $(N-1)$ dimensional balls,
the problem is treated in \cite{Pic}. In \cite{BZ3}, a continuity result is
obtained in the case $K_n \subseteq M$ and $K_n$ satisfies appropriate capacitary
conditions on the boundary. In Murat \cite{Mu}
and Del Vecchio \cite{Dv} (see also \cite{San},\cite{San2}), the case of a sieve
(Neumann sieve) is considered:
the transmission conditions that occur in the limit are determined in
relation to capacitary properties of the holes of the sieve.
\par
In this paper, we suppose that the sets $K_n$, locally, are sufficiently
regular subsets of $(N-1)$-dimensional Lipschitz submanifolds of
$\R^N$ in such a way that homogenization effects due to the possible holes
cannot occur.
\par
Let $\pi$ be the hyperplane $x_N=0$ in $\R^N$ and let $C$ be an
$(N-1)$-dimensional finite closed cone with nonempty relative interior.
We say that the sequence
$(K_n)$ satisfies the $C$-condition if there exist constants $\delta,L_1,L_2>0$
such that, for all $n$ and
for all $x \in K_n$, there exists $\Phi_x\,:\, B_{\delta}(x) \to \R^N$ with
\begin{itemize}
\item[(a)] for all $z_1,z_2 \in B_{\delta}(x)$:
$$
L_1 |z_1-z_2| \le |\Phi_x(z_1)- \Phi_x(z_2)| \le L_2 |z_1-z_2|;
$$
\item[]
\item[(b)]
$\Phi_x(x)=0$ and $\Phi_x(B_{\delta}(x) \cap K_n) \subseteq \pi$;
\item[]
\item[(c)]
for all $y \in B_{\frac{\delta}{2}}(x) \cap K_n$,
$$
\Phi_x(y) \in C_y \subseteq \Phi_x(B_{\delta}(x) \cap K_n)
$$
\end{itemize}
for some finite closed cone $C_y$ in $\pi$
congruent to $C$.
Conditions $(a)$, $(b)$ imply that, near $x$, $K_n$ is a subset of an
$(N-1)$-dimensional Lipschitz submanifold $M_{n,x}$ of $\R^N$ and condition
$(c)$ implies that $K_n$ is sufficiently regular
in $M_{n,x}$, essentially a finite union of Lipschitz subsets.
\par
The main result of the paper is that, if the sequence $(K_n)$
satisfies the $C$-condition and $K_n \to K$ in the Hausdorff metric,
then the spaces
$W^{1,p}(\Om \setminus K_n)$ converge in the sense of Mosco (see Section \ref{prel})
to the space $W^{1,p}(\Om \setminus K)$ for $1<p\le2$.
As a consequence for the case $p=2$, the problems (\ref{intpbn}) are stable, that is
transmission
conditions in the limit are avoided.
\par
The hypotheses above  are not
sufficient to cover the case $p>2$;
moreover, point $(b)$ in $C$-condition cannot be omitted:
in fact a sort of ``curvilinear'' cone condition given only by points $(a)$ and
$(c)$ does not provide the Mosco convergence. We will see these facts through
explicit examples.
\par
The paper is organized as follows: in Section \ref{prel}, we introduce the basic
notation; after some preliminaries, we prove the main
stability result in Section \ref{main}. In Section \ref{Gammaconv}, we give
the above mentioned examples of non-stability which require some basic
techniques of $\Gamma$-convergence.

\section{Notation and preliminaries}
\label{prel}

In this section, we introduce the basic notation and the tools employed in
the rest of the paper.
\vskip 10pt \noindent
{\it The Mosco convergence.}\;
Let $X$ be a reflexive Banach space, $(Y_n)$ a sequence of closed subspaces
of $X$. Let us pose
\begin{equation}
\label{liminf}
Y':= \{x \in X\,:\, x=w \hbox{-}\lim y_{n_k}\,,\,y_{n_k} \in Y_{n_k},\,
n_k \to +\infty\}
\end{equation}
and
\begin{equation}
\label{limsup}
Y'':= \{x \in X\,:\, x=s \hbox{-}\lim y_{n}\,,\,y_{n} \in Y_{n}
\; {\rm for}\;n\;{\rm large}\};
\end{equation}
$Y'$ and $Y''$ are called, respectively, the {\it weak-limsup} and the
{\it strong-liminf} of the sequence $(Y_n)$ in the sense of Mosco.
We say that the sequence
$(Y_n)$ converges in the sense of Mosco if $Y'=Y''=Y$ and we call $Y$ the
Mosco limit of $(Y_n)$. Clearly $Y'' \subseteq Y'$: as a consequence,
in order to prove
that $Y_n \to Y$ in the sense of Mosco, it is sufficient to prove that
$Y' \subseteq Y$ ({\it weak-limsup condition}) and $Y \subseteq Y''$
({\it strong-liminf condition}). Since $Y''$ is closed, the strong-liminf
condition can be established proving the inclusion $D \subseteq Y''$,
$D$ being a dense subset of $Y$.
\par
Let $\Om'$ be open and bounded in $\R^N$, $\Om_n,\Om$ open subsets of $\Om'$,
$p \in [1,+\infty]$.
We can identify the Sobolev space $W^{1,p}(\Om_n)$ with a closed subspace of
$L^p(\Om'; \R^{N+1})$ through the map
\begin{equation}
\label{identif}
\begin{array}{c}
W^{1,p}(\Om_n) \\ u
\end{array}
\begin{array}{c}
\longrightarrow \\ \longmapsto
\end{array}
\begin{array}{l}
L^p(\Om'; \R^{N+1}) \\ (u,D_1u, \dots, D_Nu)
\end{array}
\end{equation}
with the convention of extending $u$ and $\nabla u$ to zero on
$\Om' \setminus \Om_n$.
\par
Let $Y$ and $Y_n$ be the closed subspaces of $L^p(\Om'; \R^{N+1})$
corresponding to $W^{1,p}(\Om)$ and
$W^{1,p}(\Om_n)$ respectively.
We say that $W^{1,p}(\Om_n)$ converges to $W^{1,p}(\Om)$ in the sense of Mosco if
$Y$ is the Mosco limit of the sequence $(Y_n)$ in the space
$L^p(\Om'; \R^{N+1})$.

\vskip20pt\noindent
{\it Stability of Neumann problems.}\;
Let $\Om'$ be open and bounded in $\R^N$; consider the Neumann problems
\begin{equation}
\label{prelnpbs}
\left\{
\begin{array}{l}
-\Delta u_n + u_n = f \\ u \in H^1(\Om_n)
\end{array}
\right.
\end{equation}
and
\begin{equation}
\label{prelnpblim}
\left\{
\begin{array}{l}
-\Delta u + u = f \\ u \in H^1(\Om)
\end{array}
\right.
\end{equation}
with $f \in L^2(\Om')$, $\Om, \Om_n$ open subsets of $\Om'$; we intend
(\ref{prelnpbs}) and (\ref{prelnpblim}) in the usual weak sense, that is
$$
u \in H^1(\Om_n), \quad
\int_{\Om_n} \nabla u_n \nabla\varphi+
\int_{\Om_n} u\varphi=\int_{\Om_n}f\varphi
\quad\forall\varphi \in H^1(\Om_n)
$$
and
$$
u \in H^1(\Om), \quad
\int_{\Om} \nabla u \nabla\varphi +\int_{\Om} u\varphi=\int_{\Om}f\varphi
\quad\forall\varphi \in H^1(\Om).
$$
We say that the problems (\ref{prelnpbs}) converge to the problem
(\ref{prelnpblim}) if $(u_n, \nabla u_n) \to (u, \nabla u)$ strongly in
$L^2(\Om'; \R^{N+1})$ under the identification (\ref{identif}).

\vskip20pt
\noindent
{\it Hausdorff metric on compact sets.}
Let $\Om$ be open and bounded in $\R^N$.
We indicate the set of all compact subsets of $\Omb$ by $\ks(\Omb)$.
$\ks(\Omb)$ can be endowed with the
Hausdorff metric $d_H$ defined by
$$
d_H(K_1,K_2) := \max \left\{ \sup_{x \in K_1} {\rm dist}(x,K_2), \sup_{y \in
K_2} {\rm dist}(y,K_1)\right\}
$$
with the conventions ${\rm dist}(x, \emptyset)= {\rm diam}(\Om)$ and $\sup
\emptyset=0$, so that $d_H(\emptyset, K)=0$ if $K=\emptyset$ and
$d_H(\emptyset,K)={\rm diam}(\Om)$ if $K \not=\emptyset$. It turns out that
$\ks(\Omb)$ endowed with the Hausdorff metric is a compact space
(see e.g. \cite{Ro}).

\section{Some auxiliary results}
\label{aux}

In this section, we prove some results that are used in the proof of the
main theorem of the paper.
We begin recalling some properties of sets which satisfy the cone condition.
\par
Consider a closed ball $B \subseteq \R^N$ not containing $0$ and $x\in \R^N$.
The set
$$
C:=x+\{\lambda y\,:\, y \in B, 0 \le \lambda \le 1\}
$$
is called a {\it finite closed cone} in $\R^N$ with vertex at $x$.
\par
A {\it parallelepiped} with a vertex at the
origin is a set of the form
$$
P:=\left\{ \sum_{j=1}^N \lambda_jy_j \,:\, 0\le \lambda_j \le 1,\,
1 \le j \le N \right\}
$$
where $y_1, \dots, y_N$ are $N$ linearly independent vectors in $\R^N$.

\begin{definition}
\label{cone}
Let $C$ be a finite closed cone in $\R^N$ with vertex at the origin.
We say that a compact set $K \subseteq \R^N$ satisfies the cone condition
with respect to $C$ if for all $x \in K$ there exists a finite closed cone $C_x$
congruent to $C$ such that $x \in C_x \subseteq K$.
\end{definition}

If $K$ satisfies the cone condition with respect to a cone $C$, it turns out
that it is the union of the closure of a finite number of Lipschitz open
sets. In fact, the following result holds.

\begin{proposition}
\label{str}
Let $C$ be a finite closed cone in $\R^N$ with vertex at the origin and
let $K \subseteq \R^N$ be a compact set with $\diam{K} \le M$ which satisfies the
cone condition with respect to $C$.
Then for every $\rho >0$, there exist a finite number
$A_1, A_2, \cdots, A_m$ of compact subsets of $K$ with $\diam{A_j} \le \rho$
and a finite number $P_1, P_2, \cdots, P_m$ of congruent parallelepipeds with a
vertex at the origin such that:
\begin{itemize}
\item[(a)]
for all $x \in K$ there exists $1 \le i \le m$ with $P_i \subseteq C_x$;
\item[]
\item[(b)]
$K= \displaystyle \bigcup_{i=1}^m K_i$ where
$K_i= \displaystyle \bigcup_{x \in A_i} (x+P_i)$.
\end{itemize}
The number $m$ and the parallelepipeds $P_1, \dots, P_m$ depend only on $C, M, \rho$,
and not on the particular set $K$.
\par
Moreover there exists $\overline{\rho}>0$, depending only on $C$, such that for
$\rho < \overline{\rho}$, the following facts hold for all $i=1, \dots, m$:
\begin{itemize}
\item[(c)]
for every $y \in \partial K_i$, there exists $\eta>0$, an orthogonal coordinate
system $(\xi_1, \dots, \xi_n)$ and a Lipschitz function $f$ such that
$B_\eta(y) \cap K_i= B_\eta(y) \cap
\{\xi=(\xi_1,\dots,\xi_n)\,:\, \xi_n \le f(\xi_1,\dots,\xi_{n-1})\}$;
\item[]
\item[(d)]
$\tint{K_i}= \displaystyle \bigcup_{x \in A_i} \left( x+\tint{P_i} \right)$.
\end{itemize}
\end{proposition}

\begin{proof}
Properties $(a)$, $(b)$ and $(c)$ can be obtained as in the Gagliardo theorem
on the decomposition of open sets with the cone property (see \cite{Ad}, Thm. 4.8).
In particular, $\overline{\rho}$ can be chosen as the distance of the center of $P_i$
from $\partial P_i$; with this choice of $\overline{\rho}$, it turns out that,
if a ball $B$ of radius $r<\frac{\overline{\rho}}{2}$ is such that
$B \cap (x_1+P_i) \not=\emptyset$ and $B \cap (x_2+P_i) \not=\emptyset$
for some $x_1, x_2 \in A_i$, then $B$ cannot intersect relative opposite faces of
$x_1+P_i$ and $x_2+P_i$ respectively.
\par
Let us turn to the proof of point $(d)$. The inclusion
$$
\bigcup_{x \in A_i} \left( x+\tint{P_i} \right) \subseteq \tint{K_i}
$$
is immediate. Let $y \in \tint{K_i}$ and let $r<\frac{\overline{\rho}}{2}$ be
such that $B_r(y) \subseteq K_i$. There exists $x \in A_i$ such that
$y \in x+P_i$. If $y \in x+\tint{P_i}$ for some $z$, the result is obtained. Let
us suppose that $y \in x+ \partial P_i$. For every $z \in B_r(y)$, there exists
$x_z \in A_i$ with $z \in x_z+P_i$. If $y \in x_z+\tint{P_i}$, the proof is
concluded; let us assume by contradiction that $y \in x_z+\partial P_i$ for all
$z \in B_r(y)$.
Clearly $y-x_z$ cannot belong to the same face of $P_i$ as $z$ varies in $B_r(y)$
because this would contradict $z \in x_z+P_i$ for all $z \in B_r(y)$. Since
$B_r(y)$ cannot intersect
relative opposite faces of the parallelepipeds $x+P_i$ with $x \in A_i$, we
conclude that there exists a vertex $v_j$ of $P_i$ such that $y-x_z$ belongs to
a face passing through $v_j$ for all $z \in B_r(y)$. Let
$Q_j:=\{\lambda(x-v_j): x \in P_i, \lambda>0\}$ and let $y_n \to y$ be such
that $y-y_n \in \tint{Q_j}$. For $n$ large enough, since
$y \in x_{y_n}+\partial P_i$, we obtain $y_n \not \in x_{y_n}+P_i$
which is absurd. This concludes the proof of point $(d)$.
\end{proof}

Let now consider a sequence $(K_n)$ of compact subsets of $\R^N$ satisfying the
cone condition with respect to a given finite closed cone $C$ with vertex at the
origin. If $K_n$ converges to a compact set $K$ in the Hausdorff metric, clearly
$K$ satisfies the cone condition with respect to $C$.
Let $\Ps(K_n)$ be the family of all parallelepipeds contained in $K_n$ and congruent
to the parallelepipeds $P_1, \dots, P_m$ which appear in the decomposition $(b)$ of
Proposition \ref{str} and let $\Ps(K)$ be the analogous family for $K$.
Define $\Ps_r(K)$ as the subset of $\Ps(K)$ consisting of parallelepipeds $P$
such that there exists $n_k \to \infty$ and
$P^k \in \Ps(K_{n_k})$ with $P^k \to P$ in the Hausdorff metric.
Let us pose
\begin{equation}
\label{kreg}
K_r := \{ x \in K\,:\, x \in \tint{P'}, P' \in \Ps_r(K) \},
\end{equation}
and
\begin{equation}
\label{ksing}
K_s := K \setminus K_r.
\end{equation}
We call the elements of $K_r$ {\it regular points} of $K$ (relative to the
approximation given by $(K_n)$) and the elements of $K_s$ {\it singular points}
of $K$: $K_r$ is clearly an open set.

\begin{proposition}
\label{sing}
Let $C$ be a finite closed cone in $\R^N$ and
let $(K_n)$ be a sequence of compact subsets of $\R^N$ satisfying the cone
condition with respect to $C$ and converging to a compact set $K$ in the
Hausdorff metric. Then $\hs^{N-1}(K_s) < +\infty$.
\end{proposition}

\begin{proof}
Let us fix $\rho$ smaller than the constant $\overline{\rho}$ given by
Proposition \ref{str} (which does not depend on $n$).
By point $(b)$ of the
same proposition, we can write
$$
K_n = \bigcup_{i=1}^m K_n^i \quad\quad {\rm with} \quad\quad
K_n^i := \bigcup_{x \in A_n^i} (x+P_i)
$$
where $A_n^1, \dots, A_n^m$ are compact subsets of $K_n$ with
$\diam{A_n^i} \le \rho$ and $P_1, \dots, P_m$ are parallelepipeds with a vertex
at the origin. There exists $n_k \to \infty$ such that $A_{n_k}^i \to A^i$
in the Hausdorff metric for $i=1, \dots, m$: clearly $K_{n_k}^i$ converges
to $K^i:= \bigcup_{x \in A^i} (x+P^i)$ in the Hausdorff metric.
Let us prove that $\tint{K^i} \subseteq K_r$ for
$i=1, \dots, m$. Since $\diam{A^i} \le \rho$, by point $(d)$ of Proposition
\ref{str}, we have $\tint{K^i}= \bigcup_{x \in A^i} (x+\tint{P^i})$; given
$x_0 \in A^i$ and $x_{n_k} \in A_{n_k}^i$ with $x_{n_k} \to x_0$, we have that
$x_0+P^i$ is the Hausdorff limit of $x_{n_k}+P^i$. Since
$\tint{x_0+P^i}=x_0 +\tint{P^i}$, we conclude that $\tint{K^i} \subseteq K_r$ and
so $\bigcup_{i=1}^m \tint{K^i}\subseteq K_r$.
\par
By point $(c)$ of Proposition \ref{str}, we have that $K^i$
has Lipschitz boundary; we conclude that
$$
\hs^{N-1}(K_s)=
\hs^{N-1}(K \setminus K_r) \le
\sum_{i=1}^m \hs^{N-1}(\partial K^i) < +\infty.
$$
The proof is now complete.
\end{proof}

\section{The main result}
\label{main}

We now recall the main regularity assumption on the sequence $(K_n)$ of
compact subsets of $\R^N$ in order to obtain the stability result mentioned
in the Introduction. We assume $N \ge 3$.
\par
Let $\pi$ be the hyperplane $x_N=0$ in $\R^N$.

\begin{definition}
\label{ccond}
Let $C$ be a finite closed cone in $\R^{N-1}$ and let $(K_n)$ be a sequence of
compact subsets of $\R^N$. We say that $(K_n)$ satisfies the
$C$-condition
if there exist constants $\delta, L_1, L_2>0$ such that, for all $n$ and
for all $x \in K_n$, there exists $\Phi_x\,:\, B_{\delta}(x) \to \R^N$ with:
\begin{itemize}
\item[(a)] for all $z_1,z_2 \in B_{\delta}(x)$:
$$
L_1 |z_1-z_2| \le |\Phi_x(z_1)- \Phi_x(z_2)| \le L_2 |z_1-z_2|;
$$
\item[]
\item[(b)]
$\Phi_x(x)=0$ and $\Phi_x(B_{\delta}(x) \cap K_n) \subseteq \pi$;
\item[]
\item[(c)]
for all $y \in B_{\frac{\delta}{2}}(x) \cap K_n$,
$$
\Phi_x(y) \in C_y \subseteq \Phi_x(B_{\delta}(x) \cap K_n)
$$
for some finite closed cone $C_y$ in $\pi$
congruent to $C$.
\end{itemize}
\end{definition}

For technical reasons, we assume that $L_1 \diam{C} < \frac{1}{8} \delta$: this
is clearly not restrictive up to reducing $C$.
\par
We can now state the main result of the paper.

\begin{theorem}
\label{thmain}
Let $C$ be a finite closed cone in $\R^{N-1}$, $\Om$ a bounded open subset of
$\R^N$, $1<p\le 2$, $(K_n)$ a sequence of compact subsets of $\R^N$
satisfying the $C$-condition and
converging to a compact set $K$ in the Hausdorff metric. Then the spaces
$W^{1,p}(\Om \setminus K_n)$ converge to $W^{1,p}(\Om \setminus K)$
in the sense of Mosco.
\end{theorem}

In order to prove the main theorem, we need to analyze the structure of
the sets $K_n$ and $K$. This is done in the following lemmas.

\begin{lemma}
\label{deco}
Let $C$ be a finite closed cone in $\R^{N-1}$ and
let $(K_n)$ be a sequence of
compact subsets of $\R^N$ converging to $K$ in the Hausdorff metric.
Suppose that $(K_n)$ satisfies the $C$-condition. Then
there exist $m \ge 1$ such that, for $n$ large enough,
$$
K_n= \bigcup_{i=1}^m K_n^i
$$
with
$K_n^i$ compact,
$B_{\frac{\delta}{3}}(x_n^i) \cap K_n \subseteq
K_n^i \subseteq B_{\frac{\delta}{2}}(x_n^i)$ for some $x_n^i \in K_n$
such that $x_n^i \to x^i \in K$ for all $i=1, \dots, m$ and
$K \subseteq \bigcup_{i=1}^m B_{\frac{\delta}{3}}(x^i)$; moreover
$\Phi_{x_n^i}(K_n^i)$ satisfies the cone condition with respect to $C$ for all
$i=1, \dots, m$.
\end{lemma}

\begin{proof}
Since $K$ is compact, there exists a finite number of points
$x^1, \dots, x^m \in K$ such that
\begin{equation}
\label{kincl}
K \subseteq \bigcup_{i=1}^m B_{\frac{\delta}{4}}(x^i).
\end{equation}
As $K_n \to K$ in the Hausdorff metric, there exist $x_n^i \in K_n$ such that
$x_n^i \to x^i$ for $i=1, \dots, m$. For $n$ large enough, we clearly have
\begin{equation}
\label{incl}
K_n \subseteq \bigcup_{i=1}^m B_{\frac{\delta}{3}}(x_n^i).
\end{equation}
In order to conclude the proof, it is sufficient to take $K_n^i$ as the preimage
under $\Phi_{x_n^i}$ of the union of all
cones $C' \subseteq \pi$ congruent to $C$ such that
$C' \subseteq \Phi_{x_n^i}(B_{\delta}(x_n^i)\cap K_n)$ and
$C' \cap \Phi_{x_n^i}(\overline{B}_{\frac{\delta}{3}}(x_n^i)\cap K_n) \not=\emptyset$.
In fact, $K_n^i$ is compact and the inclusion
$B_{\frac{\delta}{3}}(x_n^i) \cap K_n \subseteq K_n^i$ comes directly from
the definition of $K_n^i$ and the fact that $(K_n)$ satisfies the $C$-condition;
moreover, the inclusion $K_n^i \subseteq B_{\frac{\delta}{2}}(x_n^i)$ comes from
the assumption $L_1 \diam{C} < \frac{1}{8}\delta$, and by (\ref{incl}) we have
$K_n=\bigcup_{i=1}^m K_n^i$.
Finally, by construction, $\Phi_{x_n^i}(K_n^i)$ satisfies the cone
condition with respect to $C$ for all $n$ and $i=1, \dots,m$, and by
(\ref{kincl}) we have $K \subseteq \bigcup_{i=1}^m B_{\frac{\delta}{3}}(x^i)$
which concludes the proof.
\end{proof}

\begin{lemma}
\label{decok}
Let $C$ be a finite closed cone in $\R^{N-1}$ and
let $(K_n)$ be a sequence of
compact subsets of $\R^N$ converging to $K$ in the Hausdorff metric.
Let $(K_n)$ satisfy the $C$-condition and let $K_n=\bigcup_{i=1}^m K_n^i$
according to the decomposition given by Lemma \ref{deco}.
Then, up to a subsequence, for $i=1, \dots, m$, $x_n^i \to x^i \in K$,
$K_n^i \to K^i \subseteq K$ in the Hausdorff metric,
$\Phi_{x_n^i} \to \Phi_i$
uniformly on $B_{\frac{3}{4}\delta}(x^i)$ with
\begin{itemize}
\item[]
\item[(a)] $\hfill K \subseteq \displaystyle \bigcup_{i=1}^m
B_{\frac{\delta}{3}}(x^i); \hfill$
\item[]
\item[(b)]$\hfill B_{\frac{\delta}{3}}(x^i) \cap K \subseteq K^i
\subseteq B_{\frac{3}{4}\delta}(x^i); \hfill$
\item[]
\item[(c)]$\hfill K= \displaystyle \bigcup_{i=1}^m K^i; \hfill$
\item[]
\item[(d)] for all $z_1,z_2 \in B_{\frac{3}{4}\delta}(x^i)$:
$$
L_1 |z_1-z_2| \le |\Phi_i(z_1)- \Phi_i(z_2)| \le L_2 |z_1-z_2|;
$$
\item[]
\item[(e)]
$\hfill \Phi_i(K \cap B_{\frac{3}{4}\delta}(x^i)) \subseteq \pi.\hfill$
\end{itemize}
Moreover, $\Phi_i(K^i)$ satisfies the cone condition with respect to $C$ for all
$i=1, \dots, m$.
\end{lemma}

\begin{proof}
By Lemma \ref{deco}, $x_n^i \to x^i \in K$ for all $i=1, \dots, m$ and
$K \subseteq \bigcup_{i=1}^m B_{\frac{\delta}{3}}(x^i)$; this proves point $(a)$.
Since $K_n \to K$ in the Hausdorff metric, up to a
subsequence, $K_n^i \to K^i \subseteq K$ in the Hausdorff metric for
$i=1, \dots, m$. Fix $i \in \{1, \dots, m\}$. Note that, for $n$ large enough,
$\overline{B}_{\frac{3}{4}\delta}(x^i) \subseteq B_{\delta}(x_n^i)$.
We deduce that $\Phi_{x_n^i}$ are well defined on $B_{\frac{3}{4}\delta}(x^i)$;
since they are equicontinuous and equibounded, we may assume
that $\Phi_{x_n^i} \to \Phi_i$ uniformly on $B_{\frac{3}{4}\delta}(x^i)$ with
$$
L_1 |z_1-z_2| \le |\Phi_i(z_1)- \Phi_i(z_2)| \le L_2 |z_1-z_2|
$$
for all $z_1,z_2 \in B_{\frac{3}{4}\delta}(x^i)$. This proves point $(d)$.
\par
Passing to the limit in the relations
$$
B_{\frac{\delta}{3}}(x_n^i) \cap K_n \subseteq K_n^i
\subseteq B_{\frac{\delta}{2}}(x_n^i)
$$
$$
K_n= \bigcup_{i=1}^m K_n^i
$$
$$
\Phi_{x_n^i}(K_n \cap B_{\frac{3}{4}\delta}(x_n^i)) \subseteq \pi,
$$
we obtain points $(b)$, $(c)$ and $(e)$.
\par
Finally, it is easy to see that $\Phi_i(K^i)$ satisfies the cone condition with
respect to $C$. In fact,
fix $y \in K^i$; since $K_n^i \to K^i$ in the Hausdorff metric, there exists
$y_n \in K_n^i$ with $y_n \to y$. As $\Phi_{x_n^i}(K_n^i)$ satisfies the
cone condition with respect to $C$, there exists $C_n$ finite closed cone in $\pi$
congruent to $C$ such that $\Phi_{x_n^i}(y_n) \in C_n \subseteq
\Phi_{x_n^i}(K_n^i)$. Up to a subsequence, $C_n \to C'$ in the Hausdorff metric
with $C'$ congruent to $C$. Then $\Phi_i(y) \in C' \subseteq \Phi_i(K^i)$ since
$\Phi_{x_n^i}(K_n^i) \to \Phi_i(K^i)$ in the Hausdorff metric.
\end{proof}

We can now pass to the proof of the main theorem.

\begin{proof}[\it Proof of Theorem \ref{thmain}]
Let $Y'$ and $Y''$ be the weak-limsup and the strong-liminf of the sequence
$W^{1,p}(\Om \setminus K_n)$ respectively. We have to prove that
$Y'=Y''=W^{1,p}(\Om \setminus K)$.
\par
Let us start with the inclusion
\begin{equation}
\label{wls}
Y'\subseteq W^{1,p}(\Om \setminus K).
\end{equation}
Let $(u_k)$ be a sequence in
$W^{1,p}(\Om \setminus K_{n_k})$ ($n_k \to +\infty$), and let
$v,w_1, \cdots, w_N \in L^p(\Om)$ be such that $u_k \to v$ and $D_i u_k \to w_i$
weakly in $L^p(\Om)$ for $i=1, \dots,N$ with the identification (\ref{identif}).
Since $K_{n_k} \to K$ in the Hausdorff
metric, it is readily seen that for $i=1, \dots,N$,
$w_i=D_i v$ in the sense of distributions in $\Om \setminus K$. Since
$(K_n)$ satisfies the $C$-condition, we have
$\leb^N(K)=0$; as a consequence, we get $v=0$ and $w_1, \dots, w_N=0$ a.e. on $K$,
and so we conclude that $(v,w_1, \dots, w_N)$ is the element of $L^p(\Om;\R^{N+1})$
associated to a function of $W^{1,p}(\Om \setminus K)$ according to
(\ref{identif}).
\par
We can thus pass to the inclusion
\begin{equation}
\label{sli}
W^{1,p}(\Om \setminus K) \subseteq Y'';
\end{equation}
we have to prove that, given
$u \in W^{1,p}(\Om \setminus K)$, there exists $u_n \in W^{1,p}(\Om \setminus K_n)$
such that
$(u_n, \nabla u_n) \to (u, \nabla u)$ strongly in $L^p(\Om;\R^{N+1})$.
By standard arguments on Mosco Convergence, it is sufficient to prove that,
given any
subsequence $n_j$, there exists a further subsequence $n_{j_k}$ and a
sequence
$u_k \in W^{1,p}(\Om \setminus K_{n_{j_k}})$ such that
$(u_k, \nabla u_k) \to (u, \nabla u)$ strongly in
$L^p(\Om; \R^{N+1})$.
Thus we deduce that, in order to prove (\ref{sli}), we can reason up to
subsequences.
\par
Using the decomposition given by Lemma \ref{deco}, there exists $m \ge 1$
such that
$$
K_n= \bigcup_{i=1}^m K_n^i
$$
with
$K_n^i$ compact,
$B_{\frac{\delta}{3}}(x_n^i) \cap K_n \subseteq K_n^i
\subseteq B_{\frac{\delta}{2}}(x_n^i)$ for some $x_n^i \in K_n$,
and $\Phi_{x_n^i}(K_n^i)$ satisfying the cone condition with respect
to $C$ for all $i=1, \dots, m$.
By Lemma \ref{decok}, up to a subsequence, $x_{n}^i \to x^i \in K$ for all
$i=1, \dots, m$, with $K \subseteq \bigcup_{i=1}^m B_{\frac{\delta}{3}}(x^i)$,
and $\Phi_{x_{n}^i} \to \Phi_i$ uniformly on $B_{\frac{3}{4}\delta}(x^i)$ such
that, for all $z_1,z_2 \in B_{\frac{3}{4}\delta}(x^i)$
$$
L_1 |z_1-z_2| \le |\Phi_i(z_1)-\Phi_i(z_2)| \le L_2 |z_1-z_2|.
$$
Moreover, $K_n^i \to K^i$ in the Hausdorff metric
with
$$
K=\bigcup_{i=1}^m K^i,
$$
$B_{\frac{\delta}{3}}(x^i) \cap K \subseteq K^i \subseteq
B_{\frac{3}{4}\delta}(x^i)$ and $\Phi_i(K^i)$ satisfies the cone condition
with respect to $C$ for all $i=1, \dots,m$. Finally, we have that
\begin{equation}
\label{convinpi}
\Phi_{x_{n}^i}(K_n^i) \to \Phi_i(K^i)
\end{equation}
in the Hausdorff metric for $i=1, \dots, m$.
\par
We begin proving the strong-liminf condition in the particular case in which
$u \in W^{1,p}(\Om \setminus K)$,
$\supt{u} \subset \subset B_{\frac{\delta}{3}}(x^i)$ and
\begin{equation}
\label{partcase}
\supt{u \circ \Phi_i^{-1}} \cap \pi \subseteq
[\Phi_i(K^i)]_r,
\end{equation}
where, according to (\ref{kreg}), $[\Phi_i(K^i)]_r$ denotes the set of regular
points of $\Phi_i(K^i)$ relative to the approximation (\ref{convinpi}).
Pose $w:=u \circ \Phi_i^{-1}$;
we have $w \in W^{1,p}(\Phi_i(B_{\frac{\delta}{3}}(x)) \setminus \Phi_i(K^i))$.
As in Section \ref{aux}, let $\Ps_r(\Phi_i(K^i))$ denote the family of
parallelepipeds contained in $\Phi_i(K^i)$ and congruent to the parallelepipeds
$P_1, \dots, P_m$ given by Proposition \ref{str}, that are limit in the Hausdorff
metric of parallelepipeds $P^n$ congruent to $P_1, \dots, P_m$ and contained in
$\Phi_{x_n^i}(K_n^i)$. By (\ref{kreg}) and (\ref{partcase})
there exist $D_1, \dots, D_t \in \Ps_r(\Phi_i(K^i))$ such that
$$
\supt{w} \cap \pi \subseteq \bigcup_{j=1}^t \tintpi{D_j}
$$
where $\tintpi{\cdot}$ denotes the interior relative to $\pi$.
Let $Q_j \subseteq \tintpi{D_j}$ be a
parallelepiped in $\pi$ such that $\supt{w} \cap \pi \subseteq \tintpi{Q_j}$
and let $\varepsilon>0$ be such that, posed
$U_j:= \tintpi{Q_j} \times {]-\varepsilon, \varepsilon[}$, ($j=1, \dots, t$),
$$
\bigcup_{j=1}^t U_j \subseteq \Phi_i(B_{\frac{\delta}{3}}(x^i)).
$$
Through a partition of unity associated to $\left\{ U_1, \dots, U_t, U_0\right\}$
with $U_0:= \R^N \setminus \Phi_i(K^i)$, we may write
$$
w = \sum_{j=0}^t \psi_j w,
$$
with $\psi_j \in C^{\infty}(U_j)$, $\supt{\psi_j} \subset \subset U_j$,
so that
$$
u= \sum_{j=0}^t (\psi_j \circ \Phi_i)u.
$$
Note that $\supt{(\psi_0 \circ \Phi_i)u} \cap K=\emptyset$
so that
$$
(\psi_0 \circ \Phi_i)u \in
W^{1,p}(\Om \setminus K_n)
$$
for $n$ large enough, that is $(\psi_0 \circ \Phi_i)u \in Y''$.
In order to conclude, it is thus sufficient to deal with the case
$\supt{w} \subset \subset U_j$ for $j=1, \dots, t$.
\par
Let us fix $j \in \{1, \dots, t\}$.
Set $U^+_j:=U_j \cap (\R^{N-1} \times {]0,\varepsilon[})$,
$U^-_j:=U_j \cap (\R^{N-1} \times {]-\varepsilon,0[})$, and let
$w^{\pm}:=w_{|U^{\pm}_j}$.
We have $w^{\pm} \in W^{1,p}(U_j^{\pm})$:
let $\widetilde{w}^{\pm}$ be the extension by reflection of $w^{\pm}$
on $U_j$. Note that $\supt{\widetilde{w}^{\pm}} \subset\subset U_j$.
Up to a subsequence,
$Q_j \subseteq \Phi_{x_n^i}(K_n^i)$ because $D_j \in \Ps_r(\Phi_i(K^i))$ and
$Q_j \subseteq \tintpi{D_j}$; we deduce that $U_j \setminus \Phi_{x_n^i}(K_n^i)$
has exactly two connected components that we indicate by $B^+$ and $B^-$ (note
that they do not depend on $n$ for $n$ large).
As a consequence $\Phi_{x_n^i}^{-1}(U_j) \setminus K_n$
has exactly two connected components given by $\Phi_{x_n^i}^{-1}(B^+)$ and
$\Phi_{x_n^i}^{-1}(B^-)$ respectively. Consider
$$
v_n:=
\left\{
\begin{array}{ll}
\widetilde{w}^+ \circ \Phi_i & {\rm on}\;\; \Phi_{x_n^i}^{-1}(B^+) \\ \\
\widetilde{w}^- \circ \Phi_i & {\rm on}\;\; \Phi_{x_n^i}^{-1}(B^-).
\end{array}
\right.
$$
Since $\widetilde{w}^{\pm}$ has compact support in $U_j$, we deduce that for
$n$ large enough
$$
v_n \in W^{1,p}(\Om \setminus K_n).
$$
Since $K_n^i \to K^i$ in the Hausdorff metric and
$\widetilde{w}^{\pm} \circ \Phi_i$ does not depend on $n$, $v_n \to u$ and
$\nabla v_n \to \nabla u$ a.e. in $\Om$. By the Dominated Convergence Theorem, we
deduce that $(v_n, \nabla v_n) \to (u, \nabla u)$ in $L^p(\Om; \R^{N+1})$ under
the identification (\ref{identif}). This proves $u \in Y''$ in the case $u$
satisfies (\ref{partcase}).
\par
In order to complete the proof of the theorem, we have to see that the assumption
(\ref{partcase}) is not restrictive. Consider $u \in W^{1,p}(\Om \setminus K)$.
Let $\{\varphi_1, \dots, \varphi_m, \varphi_0\}$ be a $C^{\infty}$
partition of unity associated to $B_{\frac{\delta}{3}}(x^1), \dots,
B_{\frac{\delta}{3}}(x^m), \R^N \setminus K$.
We can write
$$
u = \sum_{i=0}^m \varphi_i u.
$$
Since $\supt{\varphi_0 u} \cap K=\emptyset$, we have that
$\supt{\varphi_0 u} \cap K_n=\emptyset$ for $n$ large enough and so
$\varphi_0 u \in W^{1,p}(\Om \setminus K_n)$.
This implies $\varphi_0 u \in Y''$.
We deduce that it is not restrictive to assume
$\supt{u} \subset \subset B_{\frac{\delta}{3}}(x^i)$ for some $i=1, \dots,m$.
\par
Let us consider
$$
K_s:= \bigcup_{i=1}^m \Phi_i^{-1} \left( [\Phi_i(K^i)]_s \right)
$$
where, according to (\ref{ksing}),
$[\Phi_i(K^i)]_s$ denotes the set of singular points
of $\Phi_i(K^i)$ under the approximation (\ref{convinpi}).
By Lemma \ref{str},
we obtain
\begin{equation}
\label{crucialpoint}
\hs^{N-2}(K_s) <+\infty;
\end{equation}
by Theorem 3 in section 4.7.2 of \cite{Ev}, since $1<p\le 2$, we deduce that
$c_p(K_s, \Om)=0$, where
$$
c_p(K_s,\Om) := \inf \left\{ \int_{\Om} |\nabla u|^p \,:\, u \in W^{1,p}_0(\Om),\,
u \ge 1 \,{\rm in\;a\;neighborhood\;of\;K_s} \right\}.
$$
By standard properties of capacity, there exists a sequence
$(\psi_k)$ in $C^{\infty}_c(\R^N)$ with $\psi_k \to 0$ in $W^{1,p}(\R^N)$ and
$\psi_k \ge 1$ on a neighborhood of $K_s$. Since
$$
u= \psi_ku +(1-\psi_k)u,
$$
we deduce that the set
$$
\ds:=
\left\{ v \in W^{1,p}(\Om \setminus K): \supt{v} \cap K_s=\emptyset \right\}
$$
is dense in $W^{1,p}(\Om \setminus K) \cap L^\infty(\Om \setminus K)$ and hence in
$W^{1,p}(\Om \setminus K)$.
As observed in Section \ref{prel}, in order to prove (\ref{sli}), it is
sufficient to check the inclusion $\ds \subseteq Y''$.
If $u \in \ds$, we have that
$$
\supt{u \circ \Phi_i^{-1}} \cap \Phi_i(K^i) \subseteq [\Phi_i(K^i)]_r.
$$
Consider $V_1,V_2 \subseteq \pi$ open in the relative topology of $\pi$ and
such that
$$
\supt{u \circ \Phi_i^{-1}} \cap \Phi_i(K^i) \subset \subset V_1
\subset \subset V_2 \subset \subset [\Phi_i(K^i)]_r;
$$
let $\varepsilon>0$ with
$U_2:=V_2 \times {]-\varepsilon,\varepsilon[}
\subseteq \Phi_i(B_{\frac{\delta}{3}}(x^i))$ and set
$U_1:=V_1 \times {]-\frac{\varepsilon}{2},\frac{\varepsilon}{2}[}$.
Consider $\varphi \in C^{\infty}_c(\Phi_i^{-1}(U_2))$ with
$0 \le \varphi \le 1$ and $\varphi \equiv 1$ on $\Phi_i^{-1}(U_1)$.
Since $u \in \ds$, we deduce $\supt{(1-\varphi)u} \cap K=\emptyset$ that is
$(1-\varphi)u \in W^{1,p}(\Om \setminus K_n)$ for $n$ large enough and so
$(1-\varphi)u \in Y''$. Moreover, since
$$
\supt{(\varphi u) \circ \Phi_i^{-1}} \cap \pi \subseteq [\Phi_i(K^i)]_r,
$$
we deduce by the previous step that $\varphi u \in Y''$. We conclude
$u= \varphi u +(1-\varphi)u \in Y''$ and
the theorem is proved.
\end{proof}

From Theorem \ref{thmain} in the case $p=2$, we may deduce the stability of
the Neumann problems mentioned in the Introduction.

\begin{corollary}
\label{stability}
Let $C$ be a finite closed cone in $\R^{N-1}$,
$(K_n)$ a sequence of compact subsets of $\R^N$ satisfying the
$C$-condition and
converging to a compact set $K$ in the Hausdorff metric.
Let $\Om$ be an open and bounded subset of $\R^N$, $f \in L^2(\Om)$, and let
$u_n$ and $u$ be the solutions of the following
Neumann problems
\begin{equation}
\label{stab1}
\left\{
\begin{array}{l}
-\Delta u_n + u_n = f \\ u \in H^1(\Om \setminus K_n),
\end{array}
\right.
\end{equation}
\begin{equation}
\label{stab2}
\left\{
\begin{array}{l}
-\Delta u + u = f \\ u \in H^1(\Om \setminus K).
\end{array}
\right.
\end{equation}
Pose $u_n=0$, $\nabla u_n=0$ on $K_n \cap \Om$, and $u=0$, $\nabla u=0$ on
$K \cap \Om$.
\par
Then we have $u_n \to u$ strongly in $L^2(\Om)$ and $\nabla u_n \to \nabla u$
strongly in $L^2(\Om; \R^N)$, so that the problems (\ref{stab1}) are stable.
\end{corollary}

\begin{proof}
Let $u_n$ be the solution of (\ref{stab1}) and $u$
the solution of (\ref{stab2}). We assume the identification (\ref{identif}).
From the equation (\ref{stab1}), we have that $(u_n, \nabla u_n)$ is bounded in
$L^2(\Om; \R^{N+1})$. There exists $v \in L^2(\Om; \R^{N+1})$ such that up
to a subsequence, $(u_n,\nabla u_n) \to v$ weakly in $L^2(\Om; \R^{N+1})$.
By Theorem \ref{thmain}, we have that $H^1(\Om \setminus K_n)$ converges
to $H^1(\Om \setminus K)$ in the sense of Mosco. Thus we deduce
$v \in H^1(\Om \setminus K)$; moreover, taking
$\varphi \in H^1(\Om \setminus K)$, there exists
$\varphi_n \in H^1(\Om \setminus K_n)$ with
$(\varphi_n, \nabla \varphi_n) \to (\varphi, \nabla \varphi)$ strongly in
$L^2(\Om; \R^{N+1})$. We conclude that
\begin{eqnarray}
\label{conver}
\int_{\Om \setminus K} \nabla v \nabla \varphi
+ \int_{\Om \setminus K} v\varphi &=&
\lim_n\int_{\Om \setminus K_n} \nabla u_n \nabla \varphi_n
+ \int_{\Om \setminus K_n}u_n\varphi_n =
\\
\nonumber
&=& \lim_n \int_{\Om \setminus K_n} f \varphi_n = \\
\nonumber
&=& \int_{\Om \setminus K} f\varphi,
\end{eqnarray}
that is $v=u$. Finally, taking $\varphi_n=u_n$ and using again
(\ref{conver}), we have that
$$
||u_n||_{L^2(\Om; \R^{N+1})} \to ||u||_{L^2(\Om; \R^{N+1})}.
$$
We conclude that $(u_n,\nabla u_n) \to (u,\nabla u)$ strongly in
$L^2(\Om; \R^{N+1})$ and so the proof is complete.
\end{proof}

\begin{remark}
\label{p-laplaciano}
Similarly, under the same hypotheses of Theorem
\ref{thmain}, we can prove that the Neumann problems
\begin{equation}
\label{plaplneu}
\left\{
\begin{array}{l}
-\Delta_pu_n+|u_n|^{p-2}u_n=f \\
u_n \in W^{1,p}(\Om \setminus K_n)
\end{array}
\right.
\end{equation}
where $1<p\le2$, $\Om$ is open and bounded in $\R^N$, $f \in L^p(\Om)$ and
$\Delta_p u_n:={\rm div}(|\nabla u_n|^{p-2} \nabla u_n)$, converge to the
Neumann problem
\begin{equation}
\label{plaplneulim}
\left\{
\begin{array}{l}
-\Delta_pu+|u|^{p-2}u=f \\
u \in W^{1,p}(\Om \setminus K),
\end{array}
\right.
\end{equation}
that is $(u_n, \nabla u_n) \to (u, \nabla u)$ strongly in $L^p(\Om; \R^{N+1})$
under the identification (\ref{identif}).
\end{remark}

\begin{remark}
{\rm
The Mosco convergence proved in Theorem \ref{thmain} is the key point in order
to prove the stability of more general problems. We now breafly sketch an
application to fracture mechanics in linearly elastic bodies.
\par
For every open and bounded set $A \subseteq \R^N$, let us pose
$$
LD^{1,2}(A):=
\left\{ u \in H^1_{\rm loc}(A;\R^N)\,:\, E(u) \in L^2(A,\msim)
\right\},
$$
where $\msim$ denotes the set of symmetric matrices of order $N$ and $E(u)$
denotes the symmetric part of the gradient of $u$. Let
$|M|:= [{\rm tr}(M^2)]^{\frac{1}{2}}$ denote the standard norm in $\msim$.
\par
Consider $(K_n)$ a sequence
of compact subsets of $\R^N$ satisfying the $C$-condition with respect to a
given $(N-1)$-dimensional finite closed cone $C$ and converging to $K$ in the
Hausdorff metric.
Let $\Om$ be open and bounded in $\R^N$ and let $\partial_D \Om$ be a Lipschitz
part of $\partial \Om$. Consider $g_n,g \in H^1(\Om;\R^N)$ with $g_n \to g$
strongly and let
$$
\Gamma_n:=\left\{ u \in LD^{1,2}(\Om \setminus K_n)\,:\, u=g_n\;{\rm on}\;
\partial_D \Om \setminus K_n \right\}
$$
and
$$
\Gamma:=\left\{ u \in LD^{1,2}(\Om \setminus K)\,:\, u=g\;{\rm on}\;
\partial_D \Om \setminus K \right\}.
$$
Given the {\it Lam\'e coefficients} $\mu,\lambda$,
let $u_n \in LD^{1,2}(\Om \setminus K_n)$ be the minimum of
$$
\min_{v \in \Gamma_n} \int_{\Om \setminus K_n}
\mu|E(v)|^2+\frac{\lambda}{2}|{\rm tr}\,Ev|^2\,d\leb^N
$$
and let $u \in LD^{1,2}(\Om \setminus K)$ be the minimum of
$$
\min_{v \in \Gamma} \int_{\Om \setminus K}
\mu|E(v)|^2+\frac{\lambda}{2}|{\rm tr}\,Ev|^2\,d\leb^N.
$$
Using the Mosco convergence given by Theorem \ref{thmain} and the density
result by Chambolle \cite{Ch} (adapted to $\Om \setminus K \subseteq \R^N$),
it can be proved that $E(u_n) \to E(u)$ strongly in $L^2(\Om; \msim)$ with the
convention of considering $E(u_n)=0$ and $E(u)=0$ on $\Om \cap K_n$ and
$\Om \cap K$ respectively.
This can be interpreted as the convergence of the equilibrium deformations
for the elastic body $\Om$ with fractures $K_n$ and boundary displacements $g_n$
to the equilibrium deformation relative to the fracture $K$ and the boundary
displacement $g$.
}
\end{remark}

\section{Non-stability examples}
\label{Gammaconv}

In this section, we give two explicit examples of non-stability when the conditions
of Theorem \ref{thmain} are violated. In Example 1, we see that the $C$-condition is
not sufficient in the case $p>2$: in fact some problems related to capacity
can occur which in the case $1<p\le2$ were avoided thank to (\ref{crucialpoint}).
In Example 2, we see that a sort of uniform ``curvilinear''
cone condition for the sequence $(K_n)$ given only by points $(a)$ and $(c)$ in
the $C$-condition does not guarantee the Mosco convergence of the spaces
$W^{1,p}(\Om \setminus K_n)$ even in the case $1<p\le2$.
\vskip10pt\noindent
EXAMPLE 1. \;\;Let $Q$, $Q'$, $Q''$ be the open unit cube in $\R^N$, $\R^{N-1}$, and
$\R^{N-2}$ respectively.
For every $n \ge 1$, let us pose
$$
K_n:=
\left\{ \left[ 0, \frac{1}{2}-\frac{1}{n}\right] \cup
\left[ \frac{1}{2}+\frac{1}{n},1 \right] \right\}
\times \overline{Q''} \times \left\{\frac{1}{2}\right\}.
$$
$(K_n)$ is a sequence of compact sets in $\R^N$ whose limit
in the Hausdorff metric is
$$
K= \overline{Q'} \times \left\{\frac{1}{2}\right\}.
$$
Let us pose
$L:=\left\{\frac{1}{2}\right\} \times \overline{Q''} \times
\left\{\frac{1}{2}\right\}$,
$S_1:=Q' \times {]0,\frac{1}{2}[}$ and $S_2:=Q' \times {]\frac{1}{2},1[}$.
\par
Let $C$ be the finite closed cone in $\R^{N-1}$ determined by
$B_{\frac{1}{8}}(P)$ with $P:=(\frac{1}{8},\frac{1}{8},\dots,\frac{1}{8})$.
Clearly $(K_n)$ satisfies the $C$-condition.
\par
We claim that, if $p>2$, then the spaces $W^{1,p}(\Om \setminus K_n)$ do
not converge to $W^{1,p}(\Om \setminus K)$ in the sense of Mosco. In fact,
assuming the Mosco convergence, by Remark \ref{p-laplaciano}, we deduce that
the Neumann problems
\begin{equation}
\label{plaplexe}
\left\{
\begin{array}{l}
-\Delta_pv+|v|^{p-2}v=f \\
v \in W^{1,p}(Q \setminus K_n)
\end{array}
\right.
\end{equation}
with $f \in L^p(Q)$ converge to the problem
\begin{equation}
\label{plaplexelim}
\left\{
\begin{array}{l}
-\Delta_pv+|v|^{p-2}v=f \\
v \in W^{1,p}(Q \setminus K).
\end{array}
\right.
\end{equation}
Let $f=\chi_{S_2}$ and let $u_n$, $u$ be the
solutions of (\ref{plaplexe}) and (\ref{plaplexelim}) respectively.
We readily deduce that $u=\chi_{S_2}$; since $(u_n,\nabla u_n) \to (u,\nabla u)$
in $L^p(Q;\R^{N+1})$ under the identification (\ref{identif}),
we obtain that $u_n \to u$ strongly in $W^{1,p}(S_i)$ for $i=1,2$.
By strong convergence in $W^{1,p}(S_1)$, we get $u_n \to 0$ $c_p$-q.e. on $L$,
while from strong convergence in $W^{1,p}(S_2)$, we deduce $u_n \to 1$
$c_p$-q.e. on $L$. Since $c_p(L,Q) \not= 0$ as $p>2$, we get a contradiction: we
conclude that the Mosco convergence does not hold.
\vskip10pt\noindent
EXAMPLE 2. \;\;Let $Q$, $Q'$, $Q''$ be the open unit cube in $\R^N$, $\R^{N-1}$,
and $\R^{N-2}$ respectively. Let us write $Q=Q' \times {]0,1[}$.
For every $n \ge 1$ let us pose
$$
K_n:=
\bigcup_{i=1}^{n-1}
\left[\frac{1}{3}, \frac{2}{3} \right]
\times \overline{Q''} \times
\left\{ \frac{i}{n} \right\}.
$$
$(K_n)$ is a sequence of compact sets in $\R^N$ whose limit
in the Hausdorff metric is
$$
K= \left[\frac{1}{3}, \frac{2}{3} \right]\times \overline{Q''} \times [0,1].
$$
Let us pose $S_1:= {\left]0,\frac{1}{3}\right[} \times Q'' \times {]0,1[}$ and
$S_2:={\left]\frac{2}{3},1\right[} \times Q'' \times {]0,1[}$.
\par
Let $C$ be the finite close cone in $\R^{N-1}$ determined by
$B_{\frac{1}{6}}(P)$ with $P:=(\frac{1}{6},\frac{1}{6},\dots,\frac{1}{6})$.
Clearly there exists $\delta>0$ such that, for all $n$ and for all $x \in K_n$,
posed
$$
\Phi_x(y):= y-x,
$$
$\Phi_x: B_{\delta}(x) \to \R^N$ satisfies
conditions (a) and (c) of Definition \ref{ccond} with respect to $C$.
Observe that condition $(b)$ is not satisfied: in particular,
$\Phi_x(B_\delta(x) \cap K_n) \not \subseteq \pi$.
\par
Let $1<p\le2$ and let us consider the Neumann problems
\begin{equation}
\label{npbs}
\left\{
\begin{array}{l}
-\Delta_p v +|v|^{p-2}v =f \\
v \in W^{1,p}(Q \setminus K_n)
\end{array}
\right.
\end{equation}
with $f \in L^p(Q)$.
We claim that the problems (\ref{npbs}) do not converge to the Neumann
problem
\begin{equation}
\label{npblim}
\left\{
\begin{array}{l}
-\Delta_p v +|v|^{p-2}v =f \\
v \in W^{1,p}(Q \setminus K)
\end{array}
\right.
\end{equation}
in the sense given in Remark \ref{p-laplaciano}, that is
$(u_n,\nabla u_n) \not\to (u,\nabla u)$
strongly in $L^p(Q; \R^{N+1})$ where $u_n$ and $u$ are the solutions of problems
(\ref{npbs}) and (\ref{npblim}) respectively and the
identification (\ref{identif}) is assumed. This
implies that $W^{1,p}(Q \setminus K_n)$ does not converge to
$W^{1,p}(Q \setminus K)$ in the sense of Mosco and so it
proves that point $(b)$ in the $C$-condition cannot be omitted.
\par
We employ a $\Gamma${-}convergence technique. Let us consider the following
functionals $F_n\,:\, L^p(Q) \to [0, \infty]$ defined by
\begin{equation}
\label{funz}
F_n(z):=
\left\{
\begin{array}{ll}
\displaystyle \frac{1}{p}\int_Q |\nabla z|^p &
{\rm if}\;z \in W^{1,p}(Q \setminus K_n) \\
+\infty & {\rm otherwise.}
\end{array}
\right.
\end{equation}
We will prove that, up to a subsequence, $(F_n)$ $\Gamma$-converges with respect
to the strong topology of $L^p(Q)$ to a functional $F$ such that
if $z \in L^p(Q)$ and $F(z) < +\infty$, then
\begin{equation}
\label{reg1}
z_{|S_i} \in W^{1,p}(S_i) \quad {\rm for}\;i=1,2,
\end{equation}
\begin{equation}
\label{reg2}
z(\cdot,x_N) \in W^{1,p}(Q') \quad {\rm for}\; {\rm a.e} \; x_N \in {]0,1[}.
\end{equation}
Let us assume for the moment (\ref{reg1}) and (\ref{reg2}). Given $f \in L^p(Q)$,
the functional
$$
G(u):= \frac{1}{p}\int_Q |u|^p- \int_Q fu
$$
is a continuous perturbation of $F_n$: as a consequence,
$$
\Gamma{-}\lim_n(F_n+G)=F+G.
$$
Note that the solution $u_n$ of problem (\ref{npbs}) is precisely the
minimum of $F_n+G$: from this, we derive that for all $n$
\begin{equation}
\label{bound}
F_n(u_n)+G(u_n) \le 0.
\end{equation}
Suppose that the problems (\ref{npbs}) converge to the problem
(\ref{npblim}): then in particular, $u_n \to u$ strongly in $L^p(Q)$ where,
as usual, $u$ is extended to $0$ on $K$.
Note that $F(u)<+\infty$ because of (\ref{bound}) and the $\Gamma$-liminf
inequality. If we choose
$$
f(x):= \chi_{S_1}
$$
we conclude that $u$ is equal to $1$ on $S_1$ and equal to $0$ on $S_2$.
With the identification (\ref{identif}), we get $u=f$.
Clearly $f(\cdot,x_N) \not\in W^{1,p}(Q')$ for $x_N \in ]0,1[$ and so
we get a contradiction.
This proves that the problems (\ref{npbs}) do not converge to problem
(\ref{npblim}).
\par
In order to perform the previous argument by contradiction,
we have to prove (\ref{reg1}) and (\ref{reg2}).
This can be done in the following way. Let $z_n \to z$ strongly in $L^p(Q)$ with
\begin{equation}
\label{lsup}
F_n(z_n) \le C <+\infty.
\end{equation}
Since
$$
\frac{1}{p}\int_{S_1} |\nabla z_n|^p +\frac{1}{p}\int_{S_2} |\nabla z_n|^p \le C,
$$
we deduce that $z_{|S_i} \in W^{1,p}(S_i)$ for $i=1,2$ and so we get (\ref{reg1}).
For a.e. $x_N \in ]0,1[$, we have that $z_n(\cdot,x_N) \to z(\cdot,x_N)$
strongly in $L^p(Q')$; by (\ref{lsup}) and Fatou's lemma, we have
$$
\frac{1}{p}\int_0^1
\left( \liminf_n \int_{Q'} |\nabla z_n(y,x_N)|^p \,dy \right)\,dx_N \le C,
$$
so that for a.e. $x_N \in ]0,1[$, there exists
$C_{x_N}>0$ and a subsequence $n_k$ such that
$$
\frac{1}{p}\int_{Q'} \sum_{i=1}^{N-1}|D_iz_{n_k}(y,x_N)|^p \,dy \le C_{x_N}.
$$
We conclude that for a.e. $x_N \in [0,1]$, $z(\cdot,x_N) \in W^{1,p}(Q')$ so that
(\ref{reg2}) is proved and the proof is complete.

\bigskip
\bigskip
\centerline{ACKNOWLEDGEMENTS}
\bigskip\noindent
The author wishes to thank Gianni Dal Maso for having proposed him the problem,
and for many helpful and interesting discussions.

\end{document}